\documentclass[letter, 10pt, final, twocolumn]{article}    
\usepackage{graphicx}           %EPS Graphics
\usepackage{mathrsfs}
\usepackage{amssymb}
\usepackage{amsmath}
\usepackage{amsthm}
\usepackage{mathtools}
\usepackage{array}
\usepackage{caption}
\usepackage{subfig}
\usepackage{enumerate}
\newcommand{\be}{\begin{equation}}
\newcommand{\ee}{\end{equation}}
\everymath{\displaystyle}
\newcommand{\ba}{\begin{array}}
\newcommand{\ea}{\end{array}}

\newcommand{\bea}{\begin{equation}\begin{array}{ccl}}
\newcommand{\eea}{\end{array}\end{equation}}

\newtheorem{property}{\textit{Property}}

\newtheorem{remark}{Remark}

\title{A Survey on Output Feedback Control of Robot Manipulators with an Application to PHANToM 1.5A Haptic Device}
\author{Soulaimane Berkane% <-this % stops a space
\thanks{This work is done as a course project (UWO, ECE 9513A) under the supervision of Dr. Ilia Polushin. } 
\thanks{The author is with the Department of Electrical and Computer Engineering, University of Western Ontario, London, Ontario, Canada. {\tt\small sberkane@uwo.ca} }
}
\date{December 20, 2014}
\begin{document}

\maketitle
\thispagestyle{empty}
\pagestyle{empty}

\begin{abstract}
This paper is aimed at presenting a short survey of velocity-free control schemes for robot manipulators. A salient feature of this work is the implementation of this control schemes on a simulation model of the 3-DOF industrial PHANToM 1.5 A robot. For both the regulation and tracking problems, two classes of control algorithms without velocity measurements are considered, namely, those based on observer-controller combination and those based on approximate derivation through an $n$ order filter. The latter was shown to be more advantageous from different aspects. First the controllers based on approximate derivation are very intuitive and easy to implement in practice. Moreover, these algorithms are proved to be globally asymptotically stable, while the observer-based controllers have only local stability results. We have also considered some algorithms that are \textit{robust} to load disturbances and to unmodeled dynamics. Simulation results have confirmed the theoretical advantages and drawbacks of each algorithm. 
\end{abstract}

\section{Introduction}
Several strategies that solve the regulation and tracking control problems of robot manipulators has been extensively presented in the literature, see for instance book \cite{Spong2006} and \cite{Lewis2003}. For both the regulation and tracking problems, the main ingredient in the majority of the available controllers is the PD (Proportional-Derivative) feedback term. Generally speaking, the difference between the available algorithms is on how to compensate for the non-linearities that are present in the robot dynamics.

One drawback of these control schemes is the requirement of measurement of motor speed. Velocity measurements, generally obtained from tacho meters, are often contaminated with a considerable amount of noise. Since the noise level limits the value of the controller gains, the achievable performance is reduced.  Moreover, speed measurement increases cost, volume and weight that is why in most robotic applications, today, velocity sensors are omitted.

One way to obtain velocity feedback, without direct sensor measurements, is to use numerical differentiation of the accurate position signal. The simplicity of this \textit{ad-hoc} technique, from an implementation point of view, make it particularly useful. However, besides the fact that there is no theoretical justification to this method, this reconstruction of velocity may be inadequate for low and high speeds \cite{Belanger1992}.

An approach that has been introduced in the literature is to design an observer that reconstruct the full state of the robot (including velocity signal) from position measurement. Then, the output of the observer is fed directly to any controller that uses full state information. Even though the separation principle does not hold generally for nonlinear systems, many authors have shown that this observer-controller structure can lead at least to some local stability results. In \cite{NICOSIA90}, a nonlinear observer that reproduces the whole robot dynamics is used in a PD plus gravity compensation scheme for the set-to-point problem and also in a Lyapunov-based like tracking algorithm for robot manipulators. The authors proved that the equilibrium is locally asymptotically stable provided that the observer gain is large enough (lower bound determined by the robot parameters and the trajectory error norms). In \cite{Berghuis1993}, using passivity approach, the authors showed that adding a term proportional to the observation error to a PD plus gravity compensation controller allows us to use a \textit{linear} observer still preserving local asymptotic stability for sufficiently high gains.

Another interesting approach, with stronger stability results, was later introduced in the literature of output feedback control of nonlinear systems. This technique is known as \textit{dirty} derivation, or approximate derivation. Global asymptotic stability has been proven for the regulation problem (see the independent works of \cite{Berghuis1993a} and \cite{Kelly1993}). It is shown that by simply adding an $n$ order filter ($n$ being the number of degrees of freedom of the robot manipulator) velocity measurements are no longer needed and thus obviating the need of observers. Nevertheless, in the tracking case, only \textit{semiglobal} asymptotic stability results were obtained (see for instance \cite{Loria1995}). Since the work of Loria et al. in 1995 \cite{Loria1995}, no one was able to come up with a solution to the long standing open problem of \textit{global} asymptotic tracking without velocity feedback for robot manipulators. Many researchers in nonlinear control community were not even sure whether this problem is solvable or not, especially after the seminal article \cite{Mazenc1994}, where the concept of unboundedness observability was introduced. Surprisingly, after more than 25 years of attempts to solve the problem mentioned above, the same author Loria came back in 2013 \cite{Loria2013} to show that in fact the same controller introduced few decades ago is uniformly globally asymptotically stable. In his paper \cite{Loria2013}, the author implicitly establishes the very intuitive conjecture the damping necessary to stabilize the system may be introduced through a simple approximate-derivatives filter. Furthermore, it is proved that such a \textit{naive} control design strategy may be applied with success to systems of higher relative degree using a cascade structure of approximate differentiators.  

In this paper, we have established a short survey on output feedback control algorithms for robot manipulators for both regulation and tracking problems. The control algorithms are presented and then discussed from a theoritical point of view (structure, stability, robustness...etc). Furthermore, we have simulated, tested and compared all these output feedback control algorithms on the three-DOF industrial robot PHANToM 1.5 A haptic device from SensAble Technologies \cite{Ref2}. 

The paper is organized as follows. In section II, we give some preliminaries on robot manipulators dynamics and some of their properties. In section III, we present three output feedback regulation algorithms for robot manipulators, then we discuss the advantages and drawbacks of each one of them. Section IV is concerned about the tracking case without velocity measurements. Also three control laws are presented and discussed. In section V, we simulate and compare all the algorithms presented in this paper and conduct some tests on our PHANToM robot model. We give some concluding remarks in section VI. 
\section{Preliminaries}
The standard equations describing the dynamics of an $n$-DOF rigid robot manipulator are given by
\begin{equation}
H(q)\ddot q+C(q,\dot q)\dot q+G(q)=\tau,
\end{equation}
where $q\in\mathbb{R}^n$ is the vector of generalized joint coordinates, $H(q)$ is the positive definite $n\times n$ inertia matrix, $C(q,\dot q)\dot q\in\mathbb{R}^n$ is the vector of Coriolis and centrifugal torques, $G(q)\in\mathbb{R}^n$ is the vector of gravitational torques, and $\tau\in\mathbb{R}^n$ is the control torque input.

Let us denote the largest and smallest eigenvalues of a matrix $A\in\mathbb{R}^{n\times n}$ by $\lambda_{\mathrm{max}}^A$ and  $\lambda_{\mathrm{min}}^A$, respectively. For an $n\times 1$ vector $x$, we shall use the Euclidean norm $||x||:=\sqrt{x^\top x}$, while the norm of an $n\times n$ matrix $A$ is the corresponding induced norm $||A||:=\sqrt{\lambda_{\mathrm{max}}^{(A^\top A)}}$. The following properties can be established \cite{ArteagaPerez1998}.
\begin{property}
It holds that $\lambda_{\mathrm{min}}^H||x||^2\leq x^\top H(q)x\leq \lambda_{\mathrm{max}}^H||x||^2$, $\forall q,x\in\mathbb{R}^n,$ such that
\begin{align*}
\lambda_{\mathrm{min}}^H:=\underset{q\in\mathbb{R}^n}{\mathrm{min}}\;\lambda_{\mathrm{min}}^{H(q)}\\
\lambda_{\mathrm{max}}^H:=\underset{q\in\mathbb{R}^n}{\mathrm{max}}\;\lambda_{\mathrm{max}}^{H(q)}.
\end{align*}
\end{property}

\begin{property}
The matrix $\dot H(q)-2C(q,\dot q)$ is skew symmetric.
\end{property}
\begin{property}
There exists $0<k_c<\infty$ such that $||C(q,x)||\leq k_c||x||$ for all $x,q\in\mathbb{R}^n$.
\end{property}

\begin{property}
For all $x,y\in\mathbb{R}^n$, we have $C(q,x)y=C(q,y)x$.
\end{property}

\begin{property}
With a proper definition of the robot model
parameters, it holds that
\begin{equation}
H(q)\ddot q+C(q,\dot q)\dot q+G(q)=Y(q,\dot q,\ddot q)\theta,
\end{equation}
where $Y(q,\dot q,\ddot q)\in\mathbb{R}^{n\times p}$ is the regressor, and $\theta\in\mathbb{R}^p$ is the constant vector of parameters. 
\end{property}

\section{Output Feedback Regulation of Robot Manipulators}
\subsection{Observer-Based Design}
One approach to solve the velocity-free feedback control problem is the design of an observer that reconstructs the velocity signal from position measurements only. However, in nonlinear control theory it is well known that an observer that asymptotically reconstructs the state of a nonlinear system does not guarantee, in general, that a given stabilizing state-feedback controller will maintain its stability properties when using the estimated state instead of the true one in the feedback loop; in general a nonlinear separation principle is not valid. 

This stability problem has motivated the design of combined controller-observer schemes for robot systems using position feedback only. For instance, the authors in \cite{NICOSIA90} proposed the following observer
\begin{equation}\label{Observer}
\left\{\begin{array}{l}
\dot{\hat q}=\hat v+k_D\bar q\\
\dot{\hat v}=H^{-1}(q)\left[\tau-C(q,\dot{\hat q})\dot{\hat q}-G(q)+L\bar q\right],
\end{array}\right.
\end{equation}
where $\bar q:=q-\hat q$ denotes the observation error, $k_D$ is a positive constant and $K_P$ is a positive definite matrix. It was shown in \cite{NICOSIA90} that, for bounded velocities $||\dot q(t)||<k_q$, if the observer gain $k_D$ satisfies $$k_D>\frac{k_ck_q}{\lambda_{\mathrm{min}}^H}$$ then the equilibrium point $\bar x:=[\bar q^\top ,\dot{\bar q}^\top ]^\top =0$ of the closed loop system is asymptotically stable, and an estimated region of attraction is given by
\begin{equation*}
\mathcal{R}_o=\left\{\bar x\in\mathbb{R}^{2n}: ||\bar x||<\sqrt{\frac{\lambda_{\mathrm{min}}^L}{\lambda_{\mathrm{max}}^L}}\left(\frac{\lambda_{\mathrm{min}}^Hk_D}{k_c}-k_q\right)\right\}.
\end{equation*}

Once the observer (\ref{Observer}) is implemented, we can analyse the stability properties of the state feedback controllers, when the estimated speeds are fed instead of the true speeds of the robot manipulator. First, let us recall the seminal result by Takegaki and Arimoto \cite{Takegaki1981}, who proposed the following controller
\begin{equation}\label{Takegaki}
\tau=G(q)-K_D\dot q-K_P\tilde q
\end{equation}
where $K_D$ and $K_P$ are symmetric positive definite matrices and $\tilde q:=q-q_d$ is the trajectory error where $q_d$ is  constant. This controller consists of a gravitation compensation and a linear static state feedback which underscores its simplicity. Now, if the speed $\dot q$ is assumed not available for measurement, we can consider the following output feedback controller
\begin{equation}\label{PD+gravity+obs}
\mathrm{\bf (R1)}\;\;\;\;\tau=G(q)-K_P\tilde q-K_D\dot{\hat q}.
\end{equation}
It is shown in \cite{NICOSIA90} that if the observer gain $k_D$ satisfies
\begin{equation}\label{observer-gain}
k_D>\frac{1}{4}\frac{(\lambda_{\mathrm{max}}^{K_D})^2}{\lambda_{\mathrm{min}}^{K_D}\lambda_{\mathrm{min}}^H}
\end{equation}
then the equilibrium point $x=0$, where $x:=\left[\tilde q^\top,\dot{\tilde q}^\top,\bar q^\top,\dot{\bar q}^\top\right]^\top$, is asymptotically stable with a region of attraction given by
 \begin{multline*}
\mathcal{R}_c=\left\{ x\in\mathbb{R}^{4n}: ||x||<\frac{1}{k_c\sqrt{2}}\sqrt{\frac{\lambda_{\mathrm{min}}^{M}}{\lambda_{\mathrm{max}}^{M}}}\left(\lambda_{\mathrm{min}}^{H}k_D-\right.\right.\\\left.\left.\frac{1}{4}\frac{(\lambda_{\mathrm{max}}^{K_D})^2}{\lambda_{\mathrm{min}}^{K_D}}\right)\right\},
\end{multline*}
where $M(q)=\mathrm{bloc\;diag}\left[K_P,H(q),L_P,H(q)\right]$.

\subsection{Global Regulation via a First Order Linear Compensator}
The main drawback of the observer-based velocity-free control scheme given by (\ref{Observer})-(\ref{PD+gravity+obs}) is that the obtained convergence is local.  The authors in \cite{Berghuis1993a} proposed a different approach that allows to come up with a \textit{globally} asymptotically stable controller for the regulation problem without velocity (the same controller was independently published also in \cite{Kelly1993}).  To avoid the use of the velocity $\dot q$ in the seminal controller (\ref{Takegaki}), they have proposed the following controller
\begin{align}\label{2}\mathrm{\bf (R2)}\;\;\;\;\left\{
\begin{array}{l}
\tau=G(q)-K_D\vartheta-K_P\tilde q,\\
\vartheta=-L\int^T_0\vartheta dt+B q,
\end{array}\right.
\end{align}
where $L$ and $B$ are some symmetric positive definite matrices. This technique is called \textit{dirty derivation} since, for diagonal matrices, we can write the transfer
\begin{equation}
\vartheta=\mathrm{diag}\left\{\frac{b_is}{s+l_i}\right\}q,
\end{equation}
thus $\vartheta$ is nothing but a \textit{filtered} derivative of $q$. It can be shown that the controller (\ref{2}) achieves the global asymptotic stability of the robot system at the equilibrium $(\dot q,\tilde q,\vartheta)\equiv 0$. In this regard, consider the following Lyapunov function candidate
\begin{multline*}
V(\dot q,\tilde q,\vartheta)=\frac{1}{2}\dot q^\top H(q)\dot q+\frac{1}{2}\tilde q^\top K_P\tilde q+\frac{1}{2}\vartheta^\top K_DB^{-1}\vartheta,
\end{multline*}
whose time derivative can be shown to be
\begin{equation}\label{dV}
\dot V(\dot q,\tilde q,x)=-\vartheta^\top(K_DB^{-1}L)\vartheta\leq 0.
\end{equation}
Therefore, all signals are bounded and stable. The asymptotic stability follows from LaSalle's invariance principle. 
\subsection{Semi-Globally Asymptotically Stable Output Feedback Regulator without Gravity Compensation}
A major drawback of the above type of controllers (with gravity compensation) is the requirement of the \textit{exact} knowledge of the gravity forces, represented by the term $G(q)$. A good estimate of $G(q)$ is hardly available since the gravity force parameters depend on the payload, which is usually unknown. A mismatch in the estimation of this term leads to a shift in the equilibrium point,  and consequently to a position steady-state error \cite{Tomei1991}. High-gain feedback reduces, but does not eliminate, this error exciting on the other hand high-frequency modes and increasing the noise sensitivity.

A standard practical remedy to compensate for gravity effects is the addition of an integral term to obtain the following PID controller
\begin{equation}\label{PID}
\tau=-K_P\tilde q-K_I\int_0^T\tilde q(s)\;ds-K_D\dot q.
\end{equation} 
It was shown in \cite{Tomei1991} that this PID type controller is locally asymptotically stable when using a particular tuning of the gains. To avoid the explicit use of the velocity signal in (\ref{PID}), the authors in \cite{Ortega1995} proposed the following velocity-free scheme without gravity compensation
\begin{equation}\label{Ortega95}\mathrm{\bf (R3)}\;\;\;\;\left\{
\begin{array}{l}
\tau=-K_P\tilde q+\nu-K_D\vartheta\\
\dot\nu=-K_I(\tilde q-\vartheta),\;\;\;\nu(0)=\nu_0\in\mathbb{R}^n\\
\vartheta=\mathrm{diag}\left\{\frac{b_is}{s+l_i}\right\}q,
\end{array}\right.
\end{equation}
where $L:=\mathrm{diag}\left\{l_i\right\}, K_P, K_I$ and $K_D$ are positive definite matrices and $B:=\mathrm{diag}\left\{b_i\right\}$ with
$$
b_i>2\frac{\lambda_{\mathrm{max}}^H}{\lambda_{\mathrm{min}}^H}, \;\;\;K_P>(4k_g+1)I
$$
such that $k_g\geq \left|\left|\frac{\partial g(q)}{\partial q}\right|\right|,\;\forall q\in\mathbb{R}^n$. \\The authors in \cite{Ortega1995} showed that, given any (possibly arbitrary large) initial conditions, there exists controller gains (sufficiently small integral gain $K_I$ and sufficiently large gain $B$) that insure the asymptotic stability of the closed-loop system. This establish a semiglobal asymptotic stability of the closed-loop system, in the sense that the domain of attraction can be arbitrary enlarged with a suitable choice of the gains.

\section{Output Feedback Tracking of Robot Manipulators}

\subsection{Observer-Based Design}
As in the regulation case, once the observer (\ref{Observer}) is implemented, we can analyse the stability properties of the state feedback controllers, when the estimated speeds are fed instead of the true speeds of the robot manipulator. For instance, the authors in \cite{NICOSIA90} has considered the following output feedback tracking control law
\begin{multline}\label{tracking1}
\mathrm{\bf (T1)}\;\;\;\;\tau=H(q)\ddot q_d+C(q,\dot{\hat q})\dot q_d+G(q)\\-K_P\tilde q-K_D(\dot{\hat q}-\dot q_d).
\end{multline}
Suppose $||\dot q_d(t)||\leq k_q$, for any $t\geq 0$, then it is shown in \cite{NICOSIA90} that if the observer gain $k_D$ satisfies
$$
k_D>\frac{k_c}{\lambda_{\mathrm{min}}^H}\left[k_q+\frac{1}{4}\frac{(\lambda_{\mathrm{max}}^{K_D}+k_ck_q)^2}{\lambda_{\mathrm{min}}^{K_D}k_c}\right],
$$
then the equilibrium point $x=0$, where $x:=\left[\tilde q^\top,\dot{\tilde q}^\top,\bar q^\top,\dot{\bar q}^\top\right]^\top$, is asymptotically stable with a region of attraction given by
 \begin{multline*}
\mathcal{R}_c=\left\{ x\in\mathbb{R}^{4n}: ||x||<\frac{1}{\sqrt{2}}\sqrt{\frac{\lambda_{\mathrm{min}}^{M}}{\lambda_{\mathrm{max}}^{M}}}\left(\frac{\lambda_{\mathrm{min}}^{B}k_D}{k_c}\right.\right.\\\left.\left.-k_q-\frac{1}{4}\frac{(\lambda_{\mathrm{max}}^{K_D}+k_ck_q)^2}{\lambda_{\mathrm{min}}^{K_D}k_c}\right)\right\},
\end{multline*}
where $M(q)=\mathrm{bloc\;diag}\left[K_P,H(q),L_P,H(q)\right]$. Again, two main drawbacks of this observer-based velocity-free control scheme are: the need of an observer and the local convergence result.
\subsection{Approximate Derivation-Based Controller}
The algorithms of control without velocity measurement discussed in precedent sections amongst other (see for instance \cite{berghuis1993passivity},\cite{Canudas1990}) require the use of observers and the injection of high gains to increase the basin of attraction. Nevertheless, the authors in \cite{Loria1995} used the technique of "dirty" derivation to obviate the necessity of observers. They proposed the following tracking control algorithm

\begin{equation}\label{Loria95}\mathrm{\bf (T2)}\;\;\;\;\left\{
\begin{array}{l}
\tau=H(q)\ddot q_d+C(q,\dot q_d)\dot q_d+G(q)-K_P\tilde q\\
\hspace{5cm}-K_D\tilde\vartheta\\
\tilde\vartheta=\mathrm{diag}\left\{\frac{b_is}{s+l_i}\right\}\tilde q
\end{array}\right.
\end{equation}
where $L:=\mathrm{diag}\left\{l_i\right\}, B:=\mathrm{diag}\left\{b_i\right\}, K_D:=\mathrm{diag}\left\{k_{d_i}\right\}$ and $K_P$  are positive definite matrices such that
$$
b_i>\frac{\lambda_{\mathrm{max}}^H}{\beta\lambda_{\mathrm{min}}^H}, \;\;\;0<\beta<1.
$$
It can be shown that for bounded initial conditions, there exists \textit{always} some sufficiently large gains of the controller (\ref{Loria95}) such that we have $\lim_{t\to\infty}\tilde q(t)=0$ with a domain of attraction that can be made \textit{arbitrary} large via a suitable tuning of the gains. This result is called a \textit{semiglobal stability} and in this particular example, it is shown that for a high gain matrix $B$ the basin of attraction can be arbitrary enlarged.

During the last two decades between the period of 1990 until 2013, the global output feedback tracking of robot manipulators was a long standing open problem and a paradigm of dynamic output feedback control of nonlinear systems. Numerous attempts have been made to solve this problem and no rigorous solution was provided.
Recently, Loria \cite{Loria2013} was able to show that the same controller, given by (\ref{Loria95}), can be shown to be globally uniformly asymptotically stable provided that the gains satisfy  
$$
\frac{k_{d_m}}{2}\frac{b_m}{a_M}>k_ck_\delta,
$$
where $(.)_m$ and $(.)_M$ denote respectively, the smallest and largest elements of $(.)_i$ and $k_\delta$ is given by
$$
\mathrm{max}\left\{\underset{t\geq 0}{\mathrm{sup}} |q_d(t)|, \underset{t\geq 0}{\mathrm{sup}}  |\dot q_d(t)|, \underset{t\geq 0}{\mathrm{sup}} |\ddot q_d(t)|\right\}\leq k_\delta.
$$
This result by Loria \cite{Loria2013} solved a problem open for 25 years even if the controller was not original. It establishes that for Lagrangian systems, the fact of introducing damping through a low-pass filter, does not alter the global property of its state-feedback counter-part. 

From a practical viewpoint the controller (\ref{Loria95}) is fairly easy to implement and the sole tuning rule is both simple and practically meaningful. As a matter of fact, (\ref{Loria95}) is reminiscent of the most elementary control strategies and employs a widely-used \textit{ad hoc} alternative to a differentiator; dirty derivative. In fact, from(\ref{Loria95}) we can write the filter
$$
\tilde\vartheta=\mathrm{diag}\left\{\frac{b_i}{a_i+s}\right\}\tilde q,
$$
which is the commonly used low-pass filter to replace the unavailable derivative $\dot{\tilde q}$.

\subsection{Robust Control via Linear Estimated State Feedback}
The authors in \cite{Berghuis1994} proposed a controller, which requires only position measurements, that consists of two parts: a linear observer that generates an estimated error state (on the joint position) and a linear feedback that makes use this estimated state. Consider the linear output-feedback controller
\begin{equation}\label{Ber94}\mathrm{\bf (T3)}\;\;\;\;
\begin{array}{l}
\mathrm{{\bf Controller}}\;\left\{\tau=-K_D\dot{\hat{e}}-K_P\hat{e}\right.	\\
\mathrm{{\bf Observer}}\;\left\{\begin{array}{l}
\dot{\hat e}=w+L_D(e-\hat e)\\
\dot w=L_P(e-\hat e)
\end{array}\right.
\end{array}
\end{equation}
where $e:=q-q_d(t)$, such that $q_d(t)$ represents the desired path to be tracked by the robot manipulator and $K_P,K_D,L_P$ and $L_D$ are some symmetric positive definite matrices that represents the gains of the controller and observer. It was shown that the control scheme (\ref{Ber94}) achieves a uniform ultimate boundedness of the closed-loop system in presence of bounded \textit{load disturbances}. Moreover, by choosing the gains large enough, the trajectory error can be made arbitrary small (asymptotic stability is achieved in the limiting case of infinite gains).  

The simplicity and easiness, when implementing the above control law, is a major advantage of this linear control scheme. This is confirmed by the experimental implementation of this linear scheme that is presented in Berghuis \cite{Berghuis1994}. It is a simple linear output feedback that does not require any model knowledge. However, to counteract the ignored knowledge of the system dynamics, the stability conditions will require to apply a high gain control input torque to the robot. This is undesirable from a practical view point in terms of the energy cost. In addition, experimental results showed that the tracking performance obtained with this simple scheme is likely to be bad.
%%%%%%%%%%%%%%%%%%%%%%%%%%%%%%%%%%%%%%%%%%%%%%%%%%%%5
\section{Simulation Tests on the PHANToM 1.5A Model}
In this section, we will be simulating, testing and comparing the velocity-free control algorithms for robot manipulators in both the regulation and tracking case. In this project, we use the dynamics of the PHANToM 1.5A haptic device from SensAble Technologies \cite{Ref2}. This is a three-DOF manipulator whose dynamic equations are in the form
\begin{multline*}
\left[\begin{array}{ccc}
H_{11} &0 &0\\
0 & H_{22} &H_{23}\\
0 & H_{32} &H_{33}
\end{array}\right]\left[\begin{array}{c}
\ddot\theta_1\\
\ddot\theta_2\\
\ddot\theta_3
\end{array}\right]+\left[\begin{array}{ccc}
C_{11} &C_{12} &C_{13}\\
C_{21} & 0 &C_{23}\\
C_{31} & C_{32} &0
\end{array}\right]\\\left[\begin{array}{c}
\dot\theta_1\\
\dot\theta_2\\
\dot\theta_3
\end{array}\right]+\left[\begin{array}{c}
0\\
G_2\\
G_3
\end{array}\right]=\left[\begin{array}{c}
\tau_1\\
\tau_2\\
\tau_3
\end{array}\right],
\end{multline*}
where $H_{ij},C_{ij}$ and $G_i$'s are given in \cite{Ref2}. In order to implement the control algorithms presented in this paper we need to determine the parameters $\lambda_{\mathrm{min}}^H,\lambda_{\mathrm{max}}^H$ and $k_c$ which are inherent form the robot dynamical properties. We will also be using the same inertial parameters given in \cite{Ref2}. We have computed these parameters and we found
$$
\lambda_{\mathrm{max}}^H=0.0052,\;\;\;\lambda_{\mathrm{min}}^H=0.0003.
$$
% it can be observed from the expressions of $C_{ij}$ that
%$$
%C_{ij}=\dot q^\topa_{ij}(q),
%$$
%for some $1\times 3$ vectors $a_{ij}(q)$. Therefore,
%\begin{equation}
%\begin{split}
%||C(q,\dot q)||&=\sqrt{\sum_{i,j}C_{ij}^2}\\
%					&=\sqrt{\sum_{i,j}\dot q^\topa_{ij}(q)a_{ij}(q)^\top\dot q}\\
%					&\leq \sqrt{\lambda_\mathrm{max}^A}||\dot q||,
%\end{split}
%\end{equation}
%where $\lambda_{\mathrm{max}}^A:=\underset{q\in\mathbb{R}^n}{\mathrm{max}}\;\lambda_{\mathrm{max}}^{A(q)}$ and $A(q):=\sum_{i,j}a_{ij}(q)a_{ij}(q)^\top$. Consequently, it is clear then that $k_c=\sqrt{\lambda_{\mathrm{max}}^A}$.

We have used the 1 induced matrix norm to calculate the constant $k_c$ as follows. The definition of the 1 induced norm is
$$
||C(q,\dot q)||_1=\underset{j}{\mathrm{max}}\;\sum_i|C_{ij}|
$$
Let us define the following constants
\begin{align*}
\alpha_1&=4I_{beyy}-4I_{bezz}+ 4I_{cyy}-4I_{czz}+ 4l_1^2m_a + l_1^2m_c\\
\alpha_2&=-4I_{ayy} + 4I_{azz} − 4I_{dfyy} + 4I_{dfzz} + l_2^2m_a+ 4l
^2_3m_c\\
\alpha_3&=l_1(l_2m_a+l_3m_c).
\end{align*}
Therefore, from the expressions of $C_{ij}$, given in (\ref{2}), we can write
\begin{align*}
|C_{11}|&\leq \frac{1}{8}|\dot\theta_2|\left(2|\alpha_1|+4|\alpha_3|\right)+\frac{1}{8}|\dot\theta_3|\left(2|\alpha_2|+4|\alpha_3|\right)\\
|C_{12}|&\leq \frac{1}{8}|\dot\theta_1|\left(|\alpha_1|+4|\alpha_3|\right)\\
|C_{13}|&\leq \frac{1}{8}|\dot\theta_1|\left(|\alpha_2|+4|\alpha_3|\right)\\
|C_{21}|&\leq \frac{1}{8}|\dot\theta_1|\left(|\alpha_1|+4|\alpha_3|\right)\\
|C_{23}|&\leq \frac{1}{2}|\dot\theta_3||\alpha_3|\\
|C_{31}|&\leq \frac{1}{8}|\dot\theta_1|\left(|\alpha_2|+4|\alpha_3|\right)\\
|C_{23}|&\leq \frac{1}{2}|\dot\theta_2||\alpha_3|.
\end{align*}
Consequently, since $||\dot q||_1=|\dot\theta_1|+|\dot\theta_2|+|\dot\theta_2|\leq \sqrt{3}||\dot q||_2$ and $||C(q,\dot q)||_2\leq\sqrt{3}||C(q,\dot q)||_1$. we can verify that
\begin{equation*}
||C(q,\dot q)||_2\leq 3\left(\frac{1}{4}|\alpha_{12}|+|\alpha_3|\right)||\dot q||_2,
\end{equation*}
where $|\alpha_{12}|:=\mathrm{max}\;\{|\alpha_1|,|\alpha_2|\}$. Therefore, we can define
$$
k_c=3\left(\frac{1}{4}|\alpha_{12}|+|\alpha_3|\right)\simeq 0.0095.
$$
\subsection{Tests and Comparison of Output Feedback Regulation Algorithms}
In this section, we conduct a simulation comparison between the regulation algorithms {\bf R1, R2} and {\bf R3} given by equations (\ref{PD+gravity+obs}), (\ref{2}) and (\ref{Loria95}), respectively. The following table recap the main, theoretically derived, features of these control algorithms.
\begin{table*}
\begin{center}
{\small
\begin{tabular}{|c@{}|m{2cm}@{}|m{2cm}@{}|m{1.5cm}@{}|m{2cm}@{}|}
\hline 
 & Technique & Gravity compensation & Integral action & Asymptotic stability \\ 
\hline 
\bf R1& Observer based & Yes & No & Semi-global \\ 
\hline 
\bf R2& Dirty derivation & Yes& No& Global \\ 
\hline 
\bf R3& Dirty derivation & No & Yes & Semi-global \\ 
\hline 
\end{tabular} }
\end{center}
\caption{}
\end{table*}

We have performed our simulations, using Simulink, by considering that our PHANToM robot is required to move to the desired set point $q_d^\top=[\pi/4,\pi/2,-2\pi/3]$. The initial conditions for $q$ and $\dot q$ are assumed to be equal to zero. For the observer-based controller ({\bf R1}), the observer states are initialized at $\hat q^\top=[\pi/2, -\pi/3,0]$ and $\hat v^\top=[-2,0.5, 1]$. For controllers {\bf R2} and {\bf R3}, the virtual signals $\vartheta$ and $\nu$ are both initialized at zero. In order to conduct the performance tests on the controllers, we propose the following gains structure
 \begin{equation*}
K_P:=\alpha_PI_{3\times 3},\;\;\;\;K_D:=\alpha_DI_{3\times 3},\;\;\;\;L:=\alpha_LI_{3\times 3}.
\end{equation*}
The performance of a given algorithm is characterized via the \textit{energy} of the control input denoted $E_\tau$, which is evaluated via the integral
$$
E_\tau:=\int^T_0||\tau||^2dt,
$$
where $T$ is the simulation time. Another performance measure is the shape of the regulation error signal (overshoot, convergence rate, settling...etc). 
\subsubsection{Observer-based regulator {\bf R1}}
We start by studying the effect of the observer gain $k_D$ on the performance of controller {\bf R1}. Let $\alpha_P=5, \alpha_D=5$ and $\alpha_L=3$, the following table and figure give the effect of varying $k_D$ on the energy of the control and the response signal.
\begin{table}
\center
\begin{tabular}{|c|c|c|c|}
\hline 
$k_D$&$10^2$&$9\times10^2$ &$9\times10^3$ \\ 
\hline 
$E_\tau\times 10^{3}$&$8.4$&$83.1$ & $855$  \\ 
\hline 
\end{tabular} 
\caption{}
\end{table}

\begin{figure}
\includegraphics[width=1.2\columnwidth]{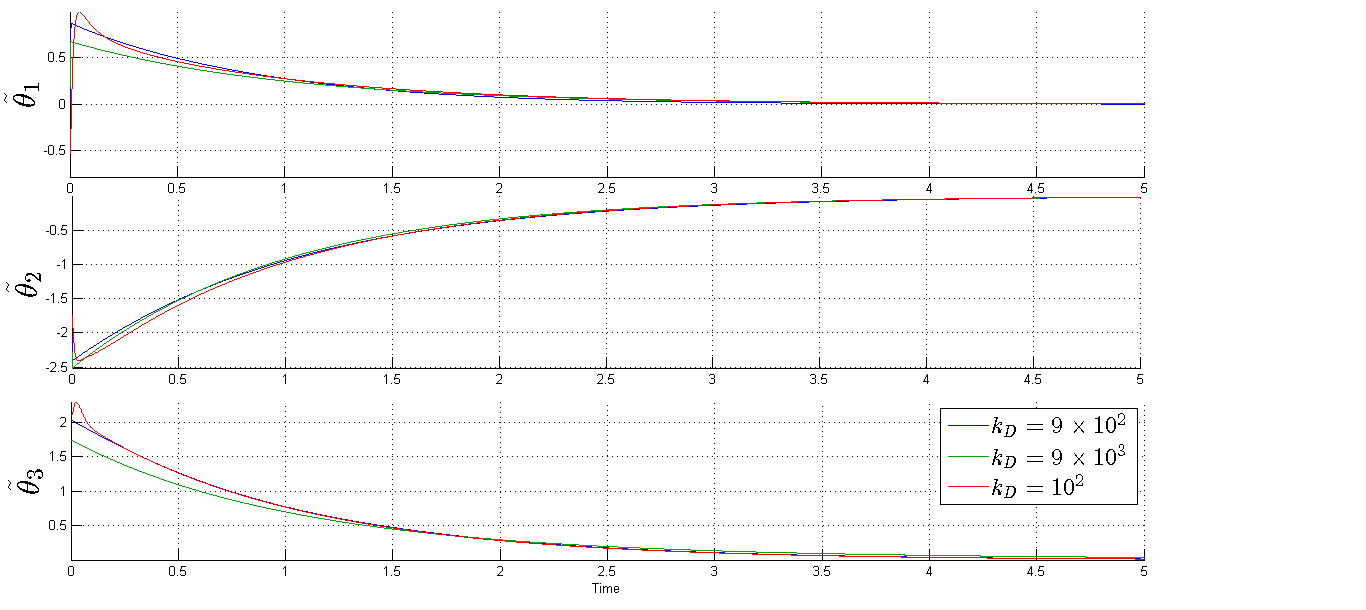}
\caption{}
\end{figure}
\begin{figure}
\includegraphics[scale=.27]{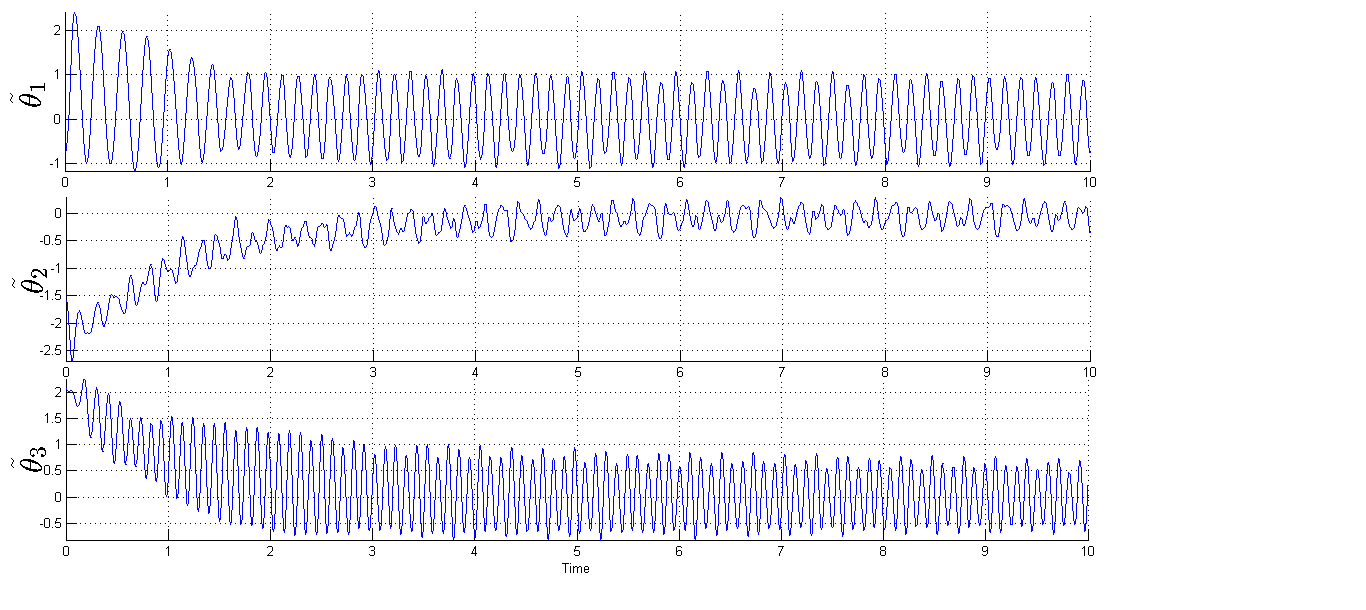}
\caption{}
\label{R1_unstable}
\end{figure}

It can be observed that the gain $k_D$ has a direct effect on the energy of the applied control. The larger $k_D$ is the more control energy we apply. Even though the value $k_D=100$ does not satisfy the stability condition (\ref{observer-gain}), our controller still be able to stabilize the system which reflects the conservative nature of inequality (\ref{observer-gain}). However, if $k_D$ is chosen very small, we have verified in simulation (see figure (\ref{R1_unstable})) that the stability is not ensured in this case, which confirm the local result of the controller {\bf R1}. 

\begin{figure}
\includegraphics[scale=.27]{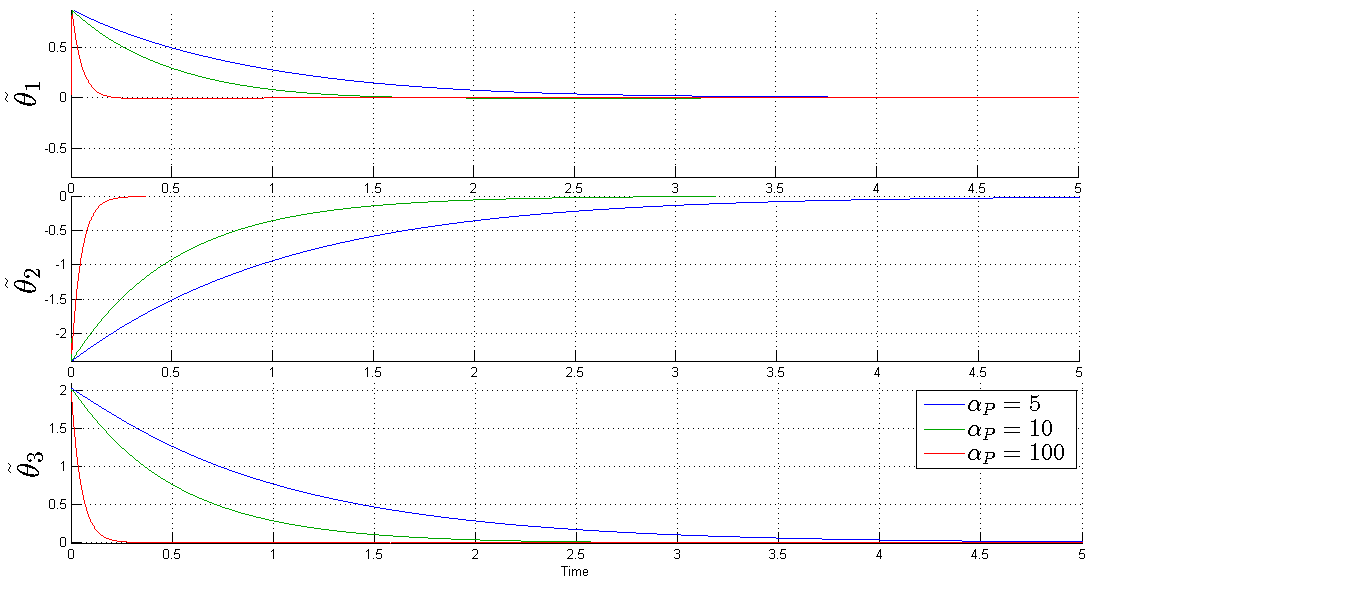}
\caption{}
\label{R1_K_P}
\end{figure}

Now let $k_D=9\times 10^2, \alpha_D=5$ and $\alpha_L=3$. Table (\ref{R1_K_P_table}) and figure (\ref{R1_K_P}) give the effect of varying $\alpha_P$ on the energy of the control and the response signal.
\begin{table}
\center
\begin{tabular}{|c|c|c|c|}
\hline 
$\alpha_P$&$5$&$10$ &$10^2$ \\ 
\hline 
$E_\tau\times 10^3$&$83.1$&$83.1$ & $84.7$  \\ 
\hline 
\end{tabular} 
\caption{}
\label{R1_K_P_table}
\end{table}

It can be seen from table (\ref{R1_K_P_table}) that the gain $K_P$ does not have a big effect on the energy of the control input. However, as seen from figure (\ref{R1_K_P}), the controller gain $K_P$ has an effect on the convergence speed. The more we increase $K_P$ the more our robot converges to the desirable configuration rapidly using almost the same amount of energy control.

\begin{figure}
\includegraphics[scale=.3]{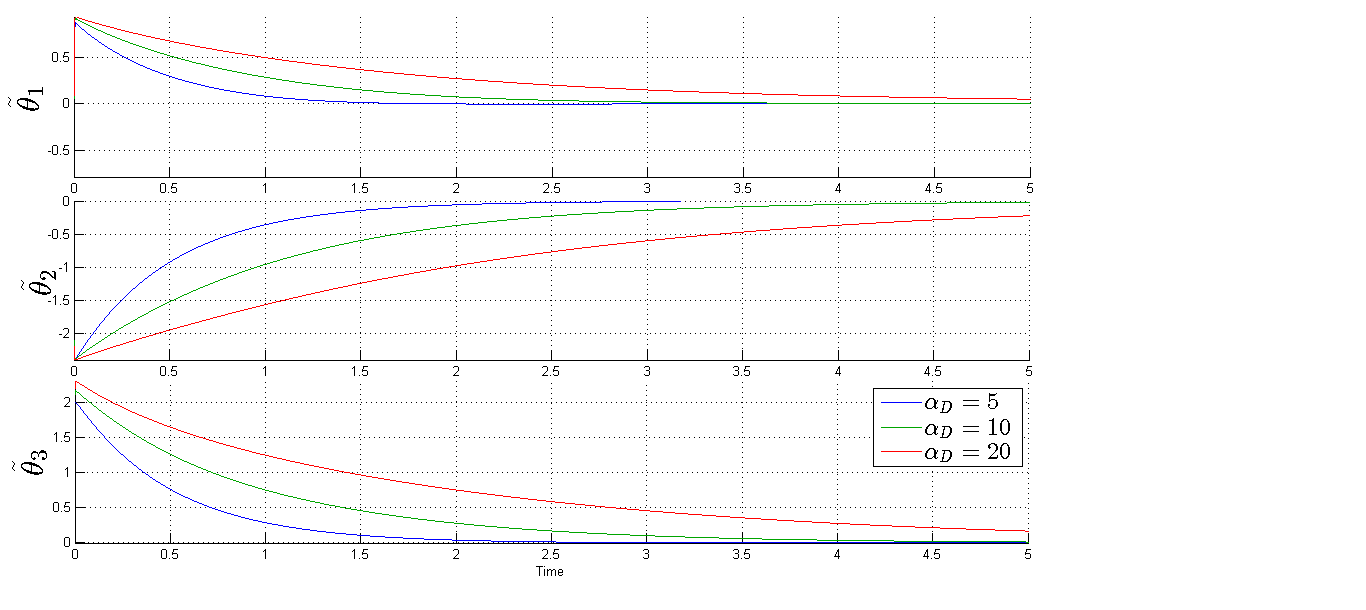}
\caption{}
\label{R1_K_D}
\end{figure}

Now let $k_D=9\times 10^2, \alpha_P=10$ and $\alpha_L=3$. Table (\ref{R1_K_D_table}) and figure (\ref{R1_K_D}) give the effect of varying $\alpha_D$ on the energy of the control and the response signal.
\begin{table}
\center
\begin{tabular}{|c|c|c|c|}
\hline 
$\alpha_D$&$5$&$10$ &$20$ \\ 
\hline 
$E_\tau\times 10^3$&$83.1$&$169.7$ & $345.3$  \\ 
\hline 
\end{tabular} 
\caption{}
\label{R1_K_D_table}
\end{table}

By increasing the gain $K_D$, the control energy is increased and the response is made more slow. Consequently, it is preferable to use a small gain $K_D$. Nevertheless, there is a minimum for the choice of the gain $\alpha_D$ that we should respect to preserve stability. In addition to this tests, we have conducted some simulations when changing the second observer gain $L$ and we have obtained no results that need to be reported. The gain $L$ is used just to ensure the convergence of the observer.  

\subsubsection{Controller {\bf R2}} 
\begin{figure}
\includegraphics[scale=.3]{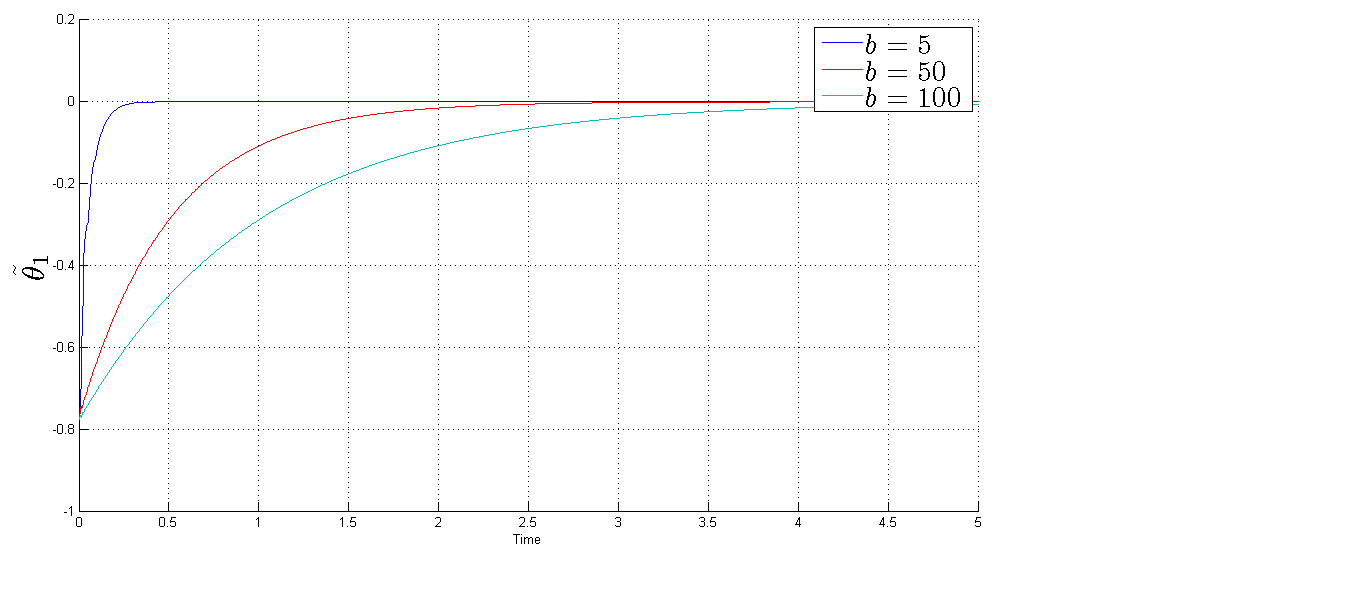}
\caption{}
\label{R2_b}
\end{figure}

Let $\alpha_P=\alpha_D=10$ and $\alpha_L=100$. We have performed multiple tests on the controller {\bf R2} with different choices of the gain $B$ and the results are given by table (\ref{R2_b_table}) and figure (\ref{R2_b}).
\begin{table}
\center
\begin{tabular}{|c|c|c|c|}
\hline 
$b$&$5$&$50$ &$100$ \\ 
\hline 
$E_\tau\times 10^3$&$.073$&$0.079$ & $0.087$  \\ 
\hline 
\end{tabular} 
\caption{}
\label{R2_b_table}
\end{table}

We can conclude from this data that the less the value of $b$ is, the faster the response is and the less the control energy is. However, it is better to not choose $b$ under a certain limit to avoid appearance of overshoots and oscillations at early times. The effect of the gain $L$ has been studied as well and the results found are totally in contrast to the effect of $B$, i.e. the larger the value of the gain $L$ the faster the response is and the less control energy is consumed. This conclusion can be intuitively derived from the time derivative of the Lyapunov function in (\ref{dV}).

\begin{figure}
\includegraphics[scale=.3]{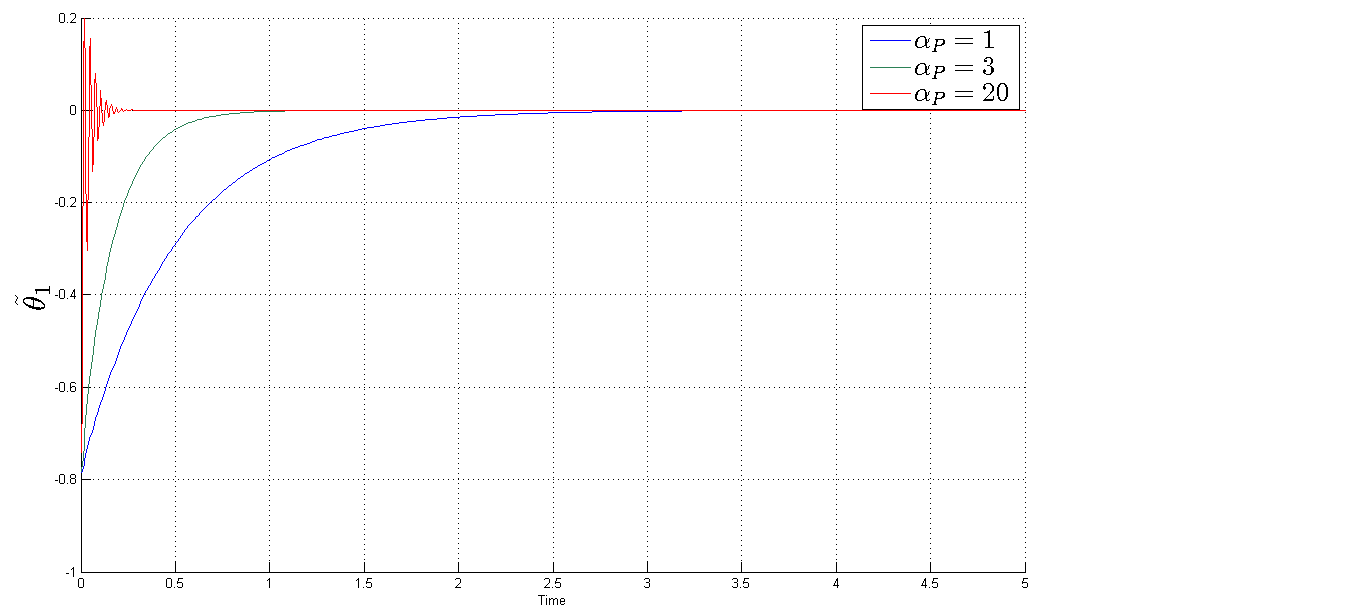}
\caption{}
\label{R2_K_P}
\end{figure}

Let $\alpha_D=10, b=5$ and $\alpha_L=100$. We have performed multiple tests on the controller {\bf R2} with different choices of the gain $K_P$ and the results are given by table (\ref{R2_K_P_table}) and figure (\ref{R2_K_P}).
\begin{table}
\center
\begin{tabular}{|c|c|c|c|}
\hline 
$\alpha_P$&$1$&$3$ &$20$ \\ 
\hline 
$E_\tau\times 10^3$&$0.006$&$0.02$ & $0.47$  \\ 
\hline 
\end{tabular} 
\caption{}
\label{R2_K_P_table}
\end{table}

As was expected, the choice of the gain $K_P$ is critical in determining the speed of convergence (rapid response for larger $K_P$). Nevertheless, we should be careful to do not exceed a certain threshold to avoid having oscillations and overshoots. This can be sensed intuitively as follows: for a larger gain $K_P$, we are ordering the controller to stabilize the robot in a very reduced and limited time and then stop suddenly at the desired reference. This cause oscillations and overshoots.

\begin{remark}
From the above simulations and tests running on both controllers {\bf R1} and {\bf R2} we noted that:
\begin{itemize}

\item For the same rise time, by comparing table (\ref{R1_K_P_table}) and table (\ref{R2_K_P_table}) it is clear that the amount of control applied by {\bf R1} is much greater than that applied by {\bf R2}.
\item Controller {\bf R1} suffers from an overshoot of more than $80\%$ (response of $\tilde\theta_1(t)$ in figure (\ref{R1_K_P})).
\item As shown in theory, controller {\bf R2} is globally asymptotically stable for any given parameter gains whereas {\bf R1} is shown to cause instability if $k_D$ is badly tuned.
\item The structure of {\bf R2} is easier from an implementation point of view than the observer-based structure of {\bf R1}.

\end{itemize}
As conclusion, controller {\bf R2} won the game challenge controller {\bf R2}.
\end{remark}
%Let $\alpha_P=3, b=5$ and $\alpha_L=100$. We have performed multiple tests on the controller {\bf R2} with different choice of the gain $K_D$ and the results are given by the following table and figure.
%\begin{table}[h!]
%\center
%\begin{tabular}{|c|c|c|c|}
%\hline 
%$\alpha_D$&$1$&$10$ &$100$ \\ 
%\hline 
%$E_\tau\times 10^3$&$0.023$&$0.020$ & $0.024$  \\ 
%\hline 
%\end{tabular} 
%\caption{}
%\end{table}
%
%\begin{figure}[h!]
%\includegraphics[scale=.3]{Simulations/R2_K_D.png}
%\caption{}
%\end{figure}
%As expected, the choice of the gain $K_P$ is critical in determining the speed of convergence (rapid response for larger $K_P$). Nevertheless, we should be careful to do not exceed a certain threshold to avoid having oscillations and overshoots. This can be sensed intuitively as follows: for larger gains $K_P$, we are ordering the controller to stabilize the robot in a very reduced time and then stop suddenly at the desired reference. This cause oscillations and overshoots.

\subsubsection{PID-like controller {\bf R3}}
\begin{figure}
\includegraphics[scale=.3]{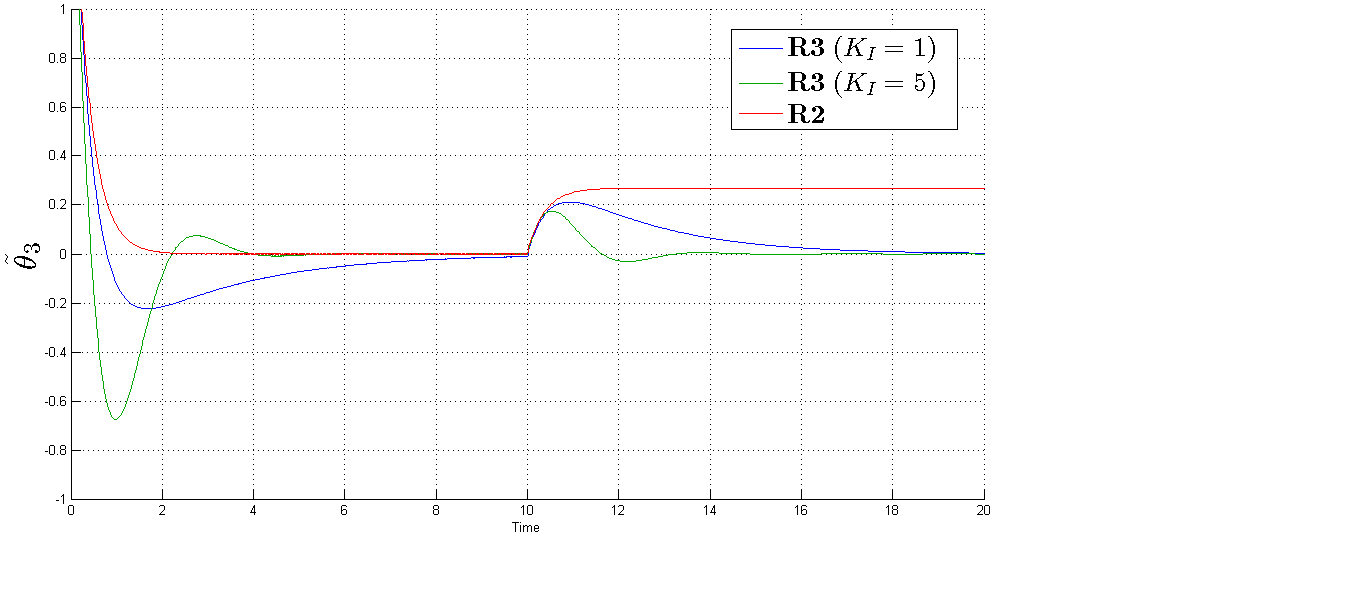}
\caption{}
\label{R23}
\end{figure}

Now we aim studying and comparing the performances of controller {\bf R2} and {\bf R3}. From a theoretical point of view, it is clear that the advantage of {\bf R3} over {\bf R2} is that the need to compensate for gravity forces is waived. In practice this is highly desirable since, for instance,  the robot can carry loads which are not necessary modelled. To illustrate this feature, we assume that our PHANToM robot carry at instant $t=5 s$ a load of mass $m_L=1 kg$ (we have added this mass to the mass $m_a$ of link A). The following figure gives a plot of the third joint error $\tilde\theta_3$ for both controllers. It can be seen from figure (\ref{R23}) that the controller {\bf R2}, even though stable, was not ale to compensate for the load disturbance applied to the robot which results in a position steady-state error. However, the PID-like regulator {\bf R3}, thanks to the integral term, was able to eliminate the effect of the disturbance load on the robot which shows the robustness of this controller.

 \subsection{Tests and Comparison of Output Feedback Tracking Algorithms}
\begin{figure}
\includegraphics[scale=.26]{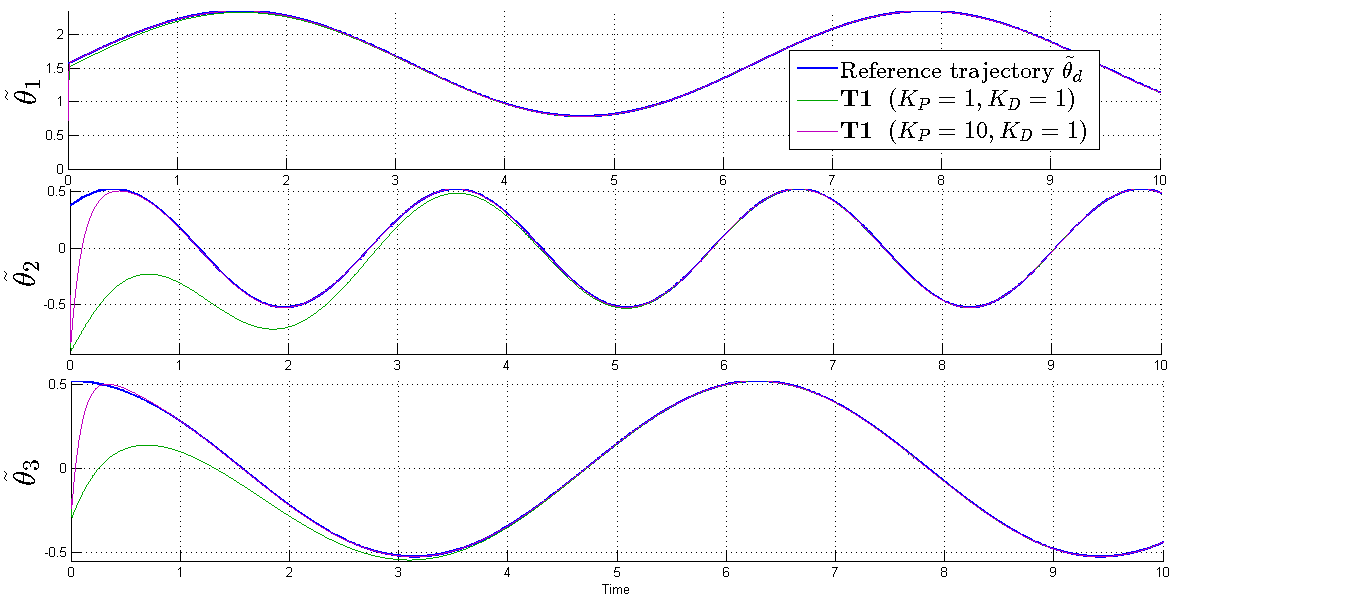}
\caption{}
\label{T1}
\end{figure}

 In this section we compare the tracking algorithms {\bf T1, T2} and {\bf T3} given by equations (\ref{tracking1}), (\ref{Loria95}) and (\ref{Loria95}), respectively. In simulating these controllers, we have considered a smooth reference trajectory given by
 $$
q_d(t)=\left[\begin{array}{c}
 \pi/4\sin(t)+\pi/2\\
 \pi/6\sin(2t+\pi/4)\\
 \pi/6\cos(t)
 \end{array}\right],
 $$
 where its first and second derivatives are given by
  $$
\dot q_d(t)=\left[\begin{array}{c}
 \pi/4\cos(t)\\
 \pi/3\cos(2t+\pi/4)\\
 -\pi/6\sin(t)
 \end{array}\right],$$
 $$\ddot q_d(t)=\left[\begin{array}{c}
-\pi/4\sin(t)\\
- 2\pi/3\sin(2t+\pi/4)\\
 -\pi/6\cos(t)
 \end{array}\right].
 $$
 The constant $k_q$ (necessary to implement controller {\bf T1} is defined by the inequality $||\dot q_d(t)||\leq k_q$. Therefore, we may chose
 $$
 k_q=\sqrt{\left(\frac{\pi}{4}\right)^2+\left(\frac{\pi}{3}\right)^2}=1.6486.
 $$
 Also the constant $k_\delta$ (necessary to implement controller {\bf T2} is defined by
 $$
\mathrm{max}\left\{\underset{t\geq 0}{\mathrm{sup}} |q_d(t)|, \underset{t\geq 0}{\mathrm{sup}}  |\dot q_d(t)|, \underset{t\geq 0}{\mathrm{sup}} |\ddot q_d(t)|\right\}\leq k_\delta.
$$
However, we have
\begin{align*}
||q_d(t)||_2\leq\sqrt{(\frac{\pi}{4}+\frac{\pi}{2})^2+(\frac{\pi}{6})^2}=2.4137\\ 
||\dot q_d(t)||_2\leq \sqrt{(\frac{\pi}{4})^2+(\frac{\pi}{3})^2}=1.6486\\
||\ddot q_d(t)||_3\leq \sqrt{(\frac{\pi}{4})^2+(2\frac{\pi}{3})^2}=2.4510.
\end{align*}
Thus $k_\delta=2.4510$.
 For the first observer-based tracking controller {\bf T1}, the tracking of the trajectory is given in figure (\ref{T1}) for different proportional gains $K_P$. For observer gain $k_D$ large enough, the controller was able to successfully track the desired trajectory after a certain time which can be set small by choosing $K_P$ large enough.

\begin{figure}
\includegraphics[scale=.29]{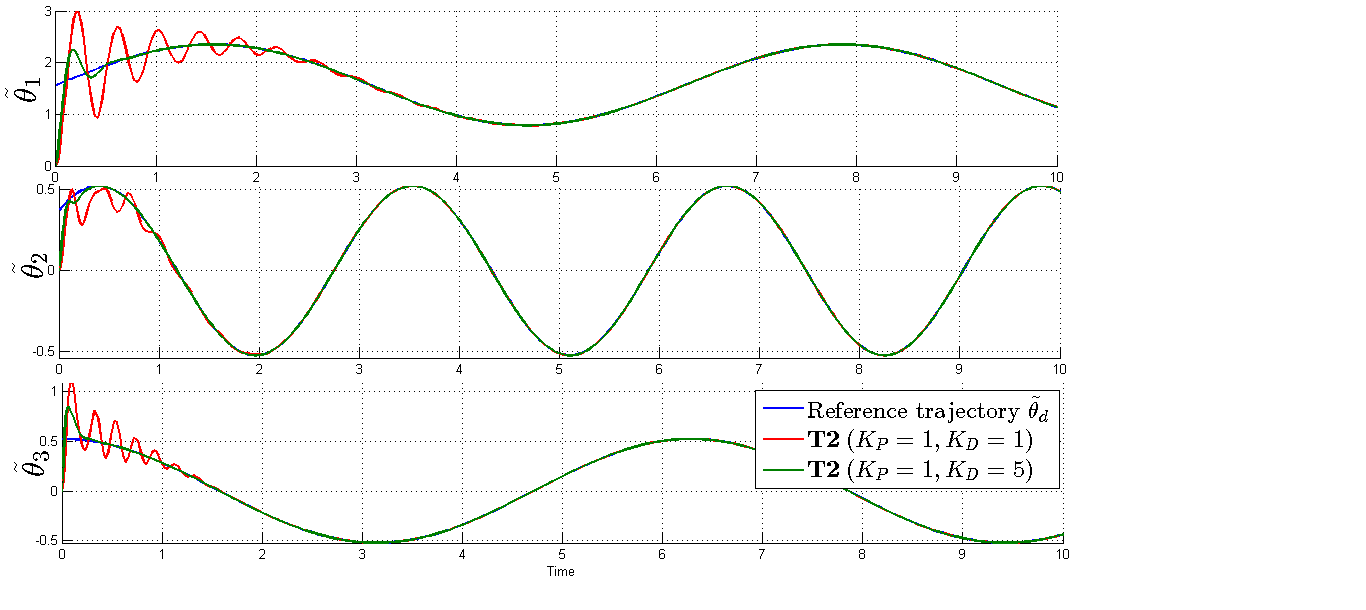}
\caption{}
\label{T2}
\end{figure}

Controller {\bf T2} was also implemented for different choices of the gain $K_D$ and the tracking trajectory is given in figure (\ref{T2}). It can be observed from this figure that the gain $K_D$ affects directly the settling time of the response. For a small gain $K_D=1$, the controlled system oscillates around the reference trajectory for a while before it converges. If we increase the gain $K_D$ to the value $5$, we can see (green response) that the convergence is more accurate and quick. Like seen in theory, the global nature of the asymptotic stability of this controller is an advantage over the observer-based controller {\bf T1}. Moreover, controller {\bf T2} is more simple form an implementation point of view.

\begin{figure}
\includegraphics[scale=.28]{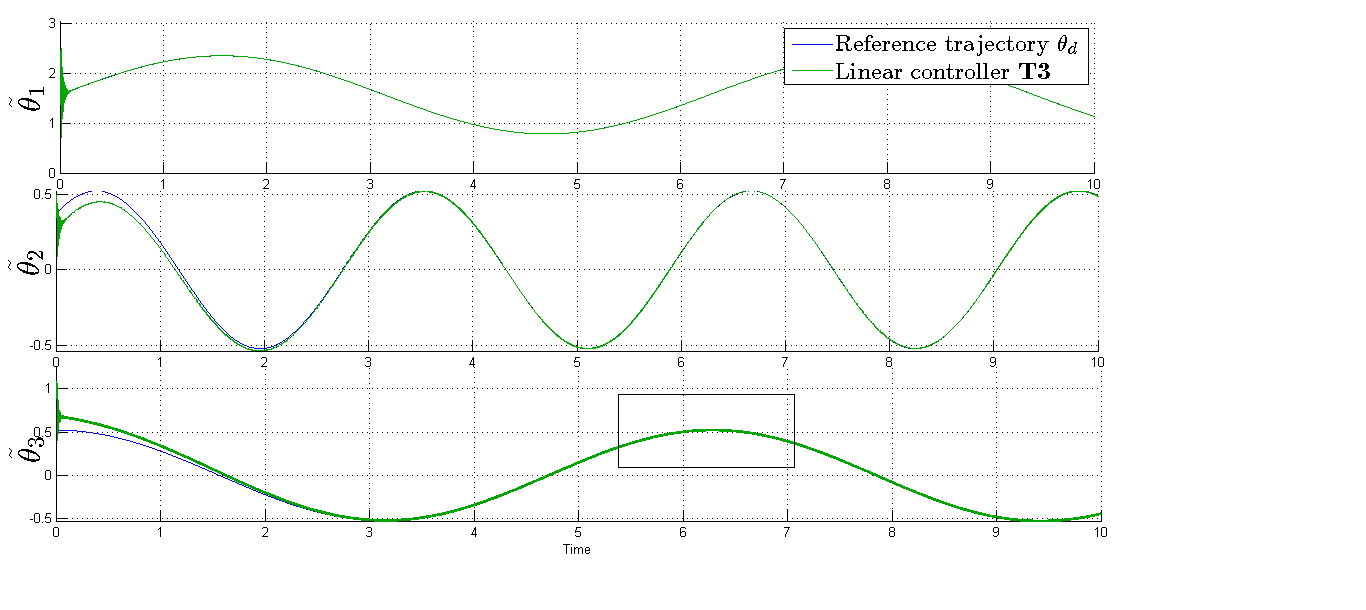}
\caption{}
\label{T3}
\end{figure}

\begin{figure}
\includegraphics[scale=.29]{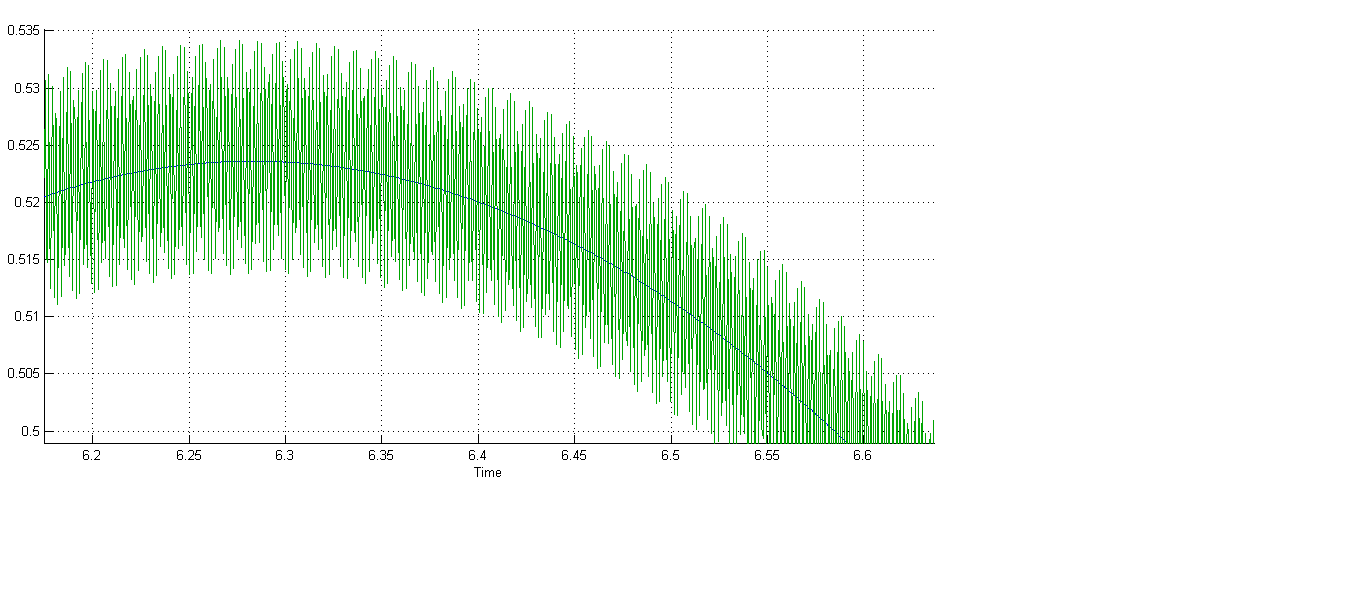}
\caption{}
\label{T3_zoom}
\end{figure}

To simulate the linear controller {\bf T3}, we have considered high gains $K_P=K_D=L_P=L_D=100$. The performance of the controller is given in figure (\ref{T3}). Apparently, the linear controller {\bf T3} was able to ensure the tracking of the reference trajectory. However, by zooming in a region of the graph, figure (\ref{T3_zoom}) shows that the controller ensures only the boundedness of the states and not the asymptotic stability (this result was proved in theory). This \textit{endless} oscillations of the states of the robot, when implementing controller {\bf T3}, is an undesirable phenomena in practice caused by neglecting the system dynamics. In fact, the motors of the robot's joints will always be chattering which may cause their damage or at least will increase the input control energy applied to the robot.  

\subsubsection{Robustness Tests}
\begin{figure}[h!]
\includegraphics[scale=.3]{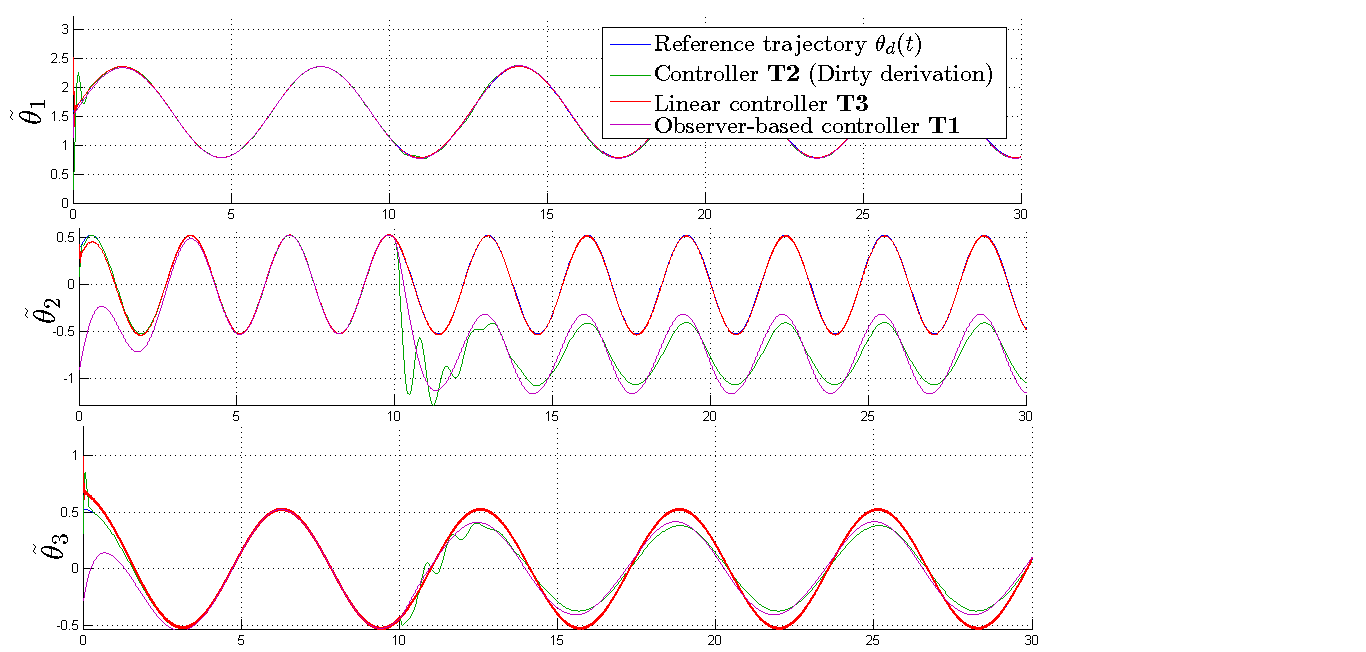}
\caption{}
\label{robustness_T}
\end{figure}

Now, we study and compare the robustness of the three proposed tracking algorithms with respect to load disturbances. We assume that at instant $t=10 s$, our robot suddenly pick up a load which will cause an increase in the value of the mass $m_a$. The result of simulation is depicted in figure (\ref{robustness_T}). It can be observed, from figure (\ref{robustness_T}), that the linear controller {\bf T3} was \textit{robust} to this load disturbance while the two other controllers {\bf T1} and {\bf T2} were not able to compensate for the additional load applied to the robot. As predicted in theory, this lack of robustness is due to fact that these controllers compensate for nonlinear terms which causes shift in the equilibrium when this terms are not exactly modelled.

\section{Conclusion}
In this paper we have presented the core algorithms of the regulation and tracking control problems without velocity measurements that are available in the literature. These algorithms were tested and compared on the model of the PHANToM 1.5A robot manipulator. Besides their local stability results, the observer-based controllers are relatively cumbersome. Moreover, we need to apply high gains to ensure our system is inside the stable region, which causes the control energy applied to the robot to be relatively large. In contrast, the controllers which are based on dirty derivation (or approximate derivation) are simple from an implementation point of view and proved to be globally asymptotically stable for both the regulation and tracking purposes. This feature give them a full credit in the area of output feedback control for nonlinear systems. Moreover, we have presented some of the robust output feedback controllers for both the regulation and tracking problems. A PID-like controller, for the regulation problem, was sufficient to provide robustness to load disturbances that are more likely to be present in a real application. For the tracking problem, we have presented a very simple output feedback \textit{linear} controller that was shown in simulation to be robust to load disturbances as well. Simulation tests were run using the PHANToM 1.5 A robot model. The comparison between these algorithms was based on two performance indexes: the energy of the control input and the shape of the system's response. 
\bibliographystyle{IEEEtran}
\bibliography{IEEEabrv,Reference}
\end{document}